\documentclass{amsart}%
\usepackage{amssymb}
\usepackage{amsfonts}
\usepackage{amsmath}
\usepackage{comment}
\usepackage{mathtools}
\usepackage{graphicx}%
\usepackage{cite}

\newtheorem{theorem}{Theorem}
\theoremstyle{plain}

\newtheorem{corollary}{Corollary}

\newtheorem{definition}{Definition}

\newtheorem{lemma}{Lemma}

\newtheorem{proposition}{Proposition}
\newtheorem{remark}{Remark}

\numberwithin{equation}{section}
\begin{document}
\title[ ]{Riccati Equation for Static Spaces and its Applications}
\author{Zhixin Wang}
\address{Department of Mathematics, Shanghai Jiao Tong University, Shanghai, 201100}
\email{jhin@sjtu.edu.cn}

\begin{abstract}
In this paper, we derive a Riccati-type equation applicable to (sub-)static Einstein spaces and examine its various applications. Specifically, within the framework of conformally compactifiable manifolds, we prove a splitting theorem for the Riemannian universal covering. Furthermore, we demonstrate two distinct methods by which the Riccati equation can establish the connectivity of the conformal boundary under the static Einstein equation. Additionally, for compact static triples possessing positive scalar curvature, we establish the compactness of the universal covering.
\end{abstract}

\maketitle

\section{Introduction}
One of the fundamental concepts in differential geometry and general relativity is that of static manifolds.
\begin{definition}
    Given a Riemannian manifold $(M,g)$ and $V\in C^{\infty}(M)$. The triple $(M,g,V)$ is called static if 
\begin{equation}
    S\coloneqq VRic+\Delta Vg-\nabla^2V=0\label{static_equ_1}
\end{equation}
And it's sub-static if
\begin{equation}
    S= VRic+\Delta Vg-\nabla^2V\geq 0\label{sub_static_equ_1}
\end{equation}
\end{definition}
\begin{comment}
The function $V$ is called potential. These spaces are interesting in both differential geometry and general relativity. For example, the negative of the left hand side of (\ref{static_equ_1}) arises as the adjoint operator of linearization of the scalar curvature, thus demonstrating its importance in prescribed scalar curvature problems. (See \cite{Co} and \cite{FM}). Secondly, they generate Riemannian manifold (possibly with singularities) or static $n+1$ space-time  of constant Ricci curvature by setting 
\begin{equation}
    h_{\pm}=\pm V^2dt^2+g\label{lifting}
\end{equation}
So its study provides with plenty of examples of $n+1$ Lorentzian or Riemannian manifold of constant Ricci curvature, see \cite{Ander}.
\end{comment}
The function $V$ is referred to as the static potential. These spaces hold significant interest in both differential geometry and general relativity. For instance, $-S$ emerges as the adjoint operator of the linearization of the scalar curvature. This reveals its critical role in problems involving prescribed scalar curvature, as discussed in references \cite{Co} and \cite{FM}.

Furthermore, these spaces can generate Riemannian manifolds (possibly with singularities) or static $n+1$ spacetimes of constant Ricci curvature by using the equation:
\begin{equation}
h_{\pm} = \pm V^2 dt^2 + g \label{lifting}
\end{equation}
Consequently, the study of such spaces provides numerous examples of 
$n+1$ Riemannian or Lorentzian manifolds with constant Ricci curvature, as noted in \cite{Ander}.

\begin{comment}Static Einstein spaces also play an important role in general relativity related problems, like positive mass theorem and no hair theorem. If a manifold is not static-Einstein, then we could perturb the metric locally, say $g_1$ so that $R_1>0$ for asymptotically flat manifold or $R_1>-n(n-1)$ in asymptotically Poincar\'e-Einstein manifold \cite{Co}. This doesn't change mass since mass is about the behavior at infinity. Then by solving a Yamabe equation we can construct a new metric $g_2$ so that $R_2=0$ or $R_2=-n(n-1)$ respectively. This procedure will decrease mass. (See \cite{Da} chapter 3 and \cite{Wo}). This argument implies that static Einstein manifold can be thought of as critical points for ADM mass or Wang mass. For explicit definition of Wang mass, readers may refer to \cite{Wa}.
\end{comment}
Static spaces also play a pivotal role in addressing problems related to general relativity, such as the positive mass theorem and the no hair theorem. For example, static spaces can be viewed as critical points for the mass. If a manifold is not static-Einstein, the metric can be locally perturbed—denoted as \( g_1 \)—to satisfy \( R_1 > 0 \) for asymptotically flat manifolds or $ R_1 > -n(n-1)$ for asymptotically Poincaré-Einstein manifolds, as discussed in \cite{Co}. Importantly, such perturbations do not alter the mass, as mass concerns the behavior at infinity. Subsequently, by solving the Yamabe equation, a new metric \( g_2 \) can be constructed such that \( R_2 = 0 \) or \( R_2 = -n(n-1) \), respectively. This process results in a reduction of mass, detailed in Chapter 3 of \cite{Da} and further discussed in \cite{Wo}. This argument underscores the interpretation of static Einstein manifolds as critical points for the ADM mass or Wang mass. For an explicit definition of Wang mass, readers are referred to \cite{Wa}.\\

\begin{comment}
One of the fundamental tool in the research of Riemannian manifold with lower bound on Ricci curvature is Riccati equation. Namely given a Riemannian manifold $(M,g)$, let $r$ be a distance function, i.e. $|\nabla r|=1$, $H$ and $A$ be the mean curvature and second fundamental form of the level sets for $r$, then
\begin{equation}
    \frac{\partial H}{\partial r}=-|A|^2-Ric(\nabla r,\nabla r)
\end{equation}
If we further assume $(M,g)$ has non-negative Ricci curvature, then together with Cauchy inequality we have $\frac{\partial H}{\partial r}\leq -\frac{1}{n-1}H^2$. This gives us the standard Ricci comparison result
\begin{equation}
    H\leq \frac{n-1}{r}\label{Ricci_comparison}
\end{equation}
\
When combined with Busemann function, we have the classical Gromoll splitting theorem \cite{CG}
\end{comment}
One of the fundamental tools in the research of Riemannian manifolds with a lower bound on Ricci curvature is the Riccati equation. Namely, given a Riemannian manifold \((M,g)\), let \(r\) be a distance function, i.e., \(|\nabla r|=1\), and let \(H\) and \(A\) be the mean curvature and second fundamental form of the level sets for \(r\), respectively. Then,
\begin{equation}
    \frac{\partial H}{\partial r} = -|A|^2 - \text{Ric}(\nabla r, \nabla r)
\end{equation}
If we further assume \((M,g)\) has non-negative Ricci curvature, for example, then together with the Cauchy inequality, we have
\begin{equation}
    \frac{\partial H}{\partial r} \leq -\frac{1}{n-1} H^2
\end{equation}
This gives us the standard Ricci comparison result:
\begin{equation}
    H \leq \frac{n-1}{r} \label{Ricci_comparison}
\end{equation}
When combined with the Busemann function, we have the classical Gromoll splitting theorem \cite{CG}:

\begin{theorem}[Cheeger-Gromoll Splitting Theorem]
Let \((M,g)\) be a complete Riemannian manifold with \( \text{Ric} \geq 0 \). If it contains a geodesic minimizing line, then it is isometric to a product
\begin{equation}
    M = \mathbb{R} \times N, \quad g = dt^2 + h
\end{equation}
where \(h\) is independent of \(t\).
\end{theorem}

Given this, one might hope to derive a Riccati-type equation for (sub-)static Einstein spaces and delve into possible applications. Such attempts have been made before, and we have the following inequality by S.Borghini and M.Fogagnolo\cite{BF}:

\begin{comment}
\begin{lemma}\label{lemma_riccati}
    Suppose $(M,,g,v)$ is a sub-static triple. Let $\tilde{g}=\frac{1}{V^2}g$, and $r$ a distance function w.r.t $\tilde{g}$, i.e. $|\tilde{\nabla} r|_{\tilde{g}}=1$, $H$ the mean curvature of level sets of $r$ computed w.r.t $g$, and set $\theta=\frac{H}{V}$. Given a $\tilde{g}$ minimizing geodesic along $\tilde{\nabla}r$, we define a new parameter $s$ satisfying 
    \begin{equation}
        ds=N^2dr\label{parameter}
    \end{equation}
    then
    \begin{equation}
        \frac{\partial \theta}{\partial s}\leq -\frac{1}{n-1}\theta^2 \label{Riccati}
    \end{equation}
\end{lemma}

In the lemma above, we construct $\tilde{g}$-distance function $r$, but compute the mean curvature of its level sets w.r.t $g$.
\end{comment}
\begin{lemma}\label{lemma_riccati}
    Suppose \((M, g, V)\) is a sub-static triple. Let \(\tilde{g} = \frac{1}{V^2} g\), and let \(\tilde{r}\) be a distance function with respect to \(\tilde{g}\), i.e., \(|\tilde{\nabla} r|_{\tilde{g}} = 1\). Let \(H\) be the mean curvature of the level sets of \(r\) computed with respect to \(g\), and set \(\theta = \frac{H}{V}\). Given a \(\tilde{g}\)-minimizing geodesic along \(\tilde{\nabla} r\), we define a new parameter \(s\) satisfying 
    \begin{equation}
        ds = V^2 dr \label{parameter}
    \end{equation}
    Then,
    \begin{equation}
        \frac{\partial \theta}{\partial s} \leq -\frac{1}{n-1} \theta^2 \label{Riccati}
    \end{equation}
    when $r$ is smooth.
\end{lemma}

In the lemma above, we start from a \(\tilde{g}\)-distance function \(r\), but compute the mean curvature of its level sets with respect to \(g\). We also perform a parameter transformation. These will be explained in section 3.

The new Riccati equation has been utilized in various scenarios. S. Brendle demonstrated that a functional is monotone along geodesic flows using this equation (\ref{Riccati}), leading to the classification of constant mean curvature surfaces in de Sitter-Schwarzschild spaces \cite{SB}. This result was later reproved using an integral formula in \cite{LX}. Furthermore, in a distinct setting for \(n=3\), G. Galloway demonstrated that no interior stable minimal 2-spheres exist, which ensures the connectivity of the external boundary \cite{Ga}. For a similar discussion, refer to \cite{Am}. This paper will delve deeper into its compelling applications.\\

It is well known that the scalar curvature of a static Einstein manifold is constant. If the scalar curvature is negative, the primary objects of interest are conformally compactifiable manifolds, which we will define in the subsequent section. For these manifolds, we have the following splitting theorem attributed to G. Galloway, S. Surya, and E. Woolgar \cite{GSW}:
\begin{theorem}[G.Galloway, S.Surya, E.Woolgar]\label{thm_splitting}
    Consider a conformally compactifiable manifold $(M^n,g)$ with conformal boundary $\Sigma$. Suppose the static Einstein equation (\ref{static_equ_1}) holds for $V$, and $\tilde{V}=\frac{1}{V}$ is a boundary defining function. We also assume \textbf{condition C} hold: the second fundamental form $II$ for level sets $\{\tilde{V}=c\}$ is positive semi-definite for small $c$. Then the Riemannian universal cover $(\tilde{M}^*,\tilde{g}^*)$ of $(\tilde{M},\tilde{g}=\frac{1}{V^2}g)$ splits isometrically as
    \begin{equation}
        \tilde{M}^*=\mathbb{R}^k\times\Sigma,\quad \tilde{g}^*=g_{\mathbb{E}}+ \tilde{h}
    \end{equation}
    where $(\mathbb{R}^k,g_{\mathbb{E}})$ is standard $k$-dimensional Euclidean space with $0\leq k\leq n$, and $(\Sigma,\tilde{h})$ is a compact Riemannian manifold with non-empty boundary. Furthermore, both $h$ and $V^*$ are independent of $r$. As a result, the Riemannian universal cover $(M^*,g^*)$ of $(M,g)$ splits isometrically as
    \begin{equation}
        M^*=\mathbb{R}^k\times \Sigma,\quad g^*=(V^{*2}g_{\mathbb{E}})+ h
    \end{equation}
    where $V^*=V\circ \pi$($\pi=$ covering map) is constant along $\mathbb{R}^k$.
\end{theorem}
The original proof in \cite{GSW} utilizes techniques from Lorentzian geometry, lifting the Riemannian manifold $(M^n, g)$ to a $(n+1)$-dimensional spacetime as specified in (\ref{lifting}). If the universal cover $(M^*,g^*)$ is not compact, constructing a line in $(M, g)$, which is lifted to a null geodesic, allows for the application of the null splitting theorem from \cite{Ga1}, achieving the desired outcome.

We propose an alternative approach by integrating (\ref{Riccati}) with the Busemann function, enabling us to rederive the result without employing a Lorentz metric. Intriguingly, this methodology closely mirrors the renowned Cheeger-Gromoll splitting theorem for manifolds with nonnegative Ricci curvature, thereby highlighting notable similarities. Additionally, by employing a small technique, we can extend the majority of this result to sub-static spaces.\\

The Riccati equation also aids in studying the topology of conformally compactifiable sub-static Einstein manifolds \cite{BF}.
\begin{theorem}[S.Borghini, M.Fogagnolo]\label{thm_connect}
    Let $(M, g, V)$ be a conformally compactifiable manifold satisfying the sub-static inequality (\ref{sub_static_equ_1}), and suppose $\tilde{V} = \frac{1}{V}$ serve as a boundary defining function. Then, the conformal boundary of $(M, g, V)$ is connected.
\end{theorem}
As previously noted, (\ref{Riccati}) has been utilized to study interior stable minimal surfaces, yielding results similar to those reported in \cite{Am} and \cite{Ga}. The recent proof introduced in \cite{BF} shifts focus to geodesics, offering a simpler approach. In this paper, we will delve into how \(\tilde{g}\) geodesic flows can be employed to tackle this issue in both of these contexts. For the convenience of readers, I will present a slightly different, yet essentially equivalent, proof to that found in \cite{BF}.

It is noteworthy that for a conformally compactifiable manifold with a specified lower Ricci bound, the following result is documented in \cite{CG2} and \cite{WY}:
\begin{theorem}[M.Cai, G.Galloway; E.Witten, S.T.Yau]\label{Ricci_connectedness}
    Let $(M, g)$ be a conformally compactifiable manifold with $Ric \geq -(n-1)g$. If one connected component of its conformal boundary possesses a non-negative Yamabe invariant, then $H_{n-1}(M, \mathbb{Z}) = 0$. Consequently, the conformal boundary of $M$ is connected.
\end{theorem}
\begin{comment}
The proof of E.Witten and Yau in \cite{WY} used a brane action $L(\Sigma)=Area(\Sigma)-Vol(\Sigma)$ where $Vol(\Sigma)$ is the volume enclosed by $\Sigma$. They showed that there exists a minimizer for $L(\Sigma)$ in each homology class provided that one connected boundary component has positive scalar curvature, and a contradiction follows. This proof is similar, in certain sense, to the proof in \cite{Ga}. While the proof in \cite{CG} used Riccati equation, which is similar to the proof in this paper.
\end{comment}
The proof by E. Witten and Yau in \cite{WY} employs a brane action $L(\Sigma) = \text{Area}(\Sigma) - \text{Vol}(\Sigma)$, where $\text{Vol}(\Sigma)$ represents the volume enclosed by $\Sigma$. They demonstrated that a minimizer for $L(\Sigma)$ exists in each homology class when one connected boundary component exhibits positive scalar curvature, leading to a subsequent contradiction. This approach bears a resemblance, in certain respects, to the method used in \cite{Ga}. Conversely, the proof in \cite{CG} utilized the Riccati equation, aligning more closely with the methodology presented in this paper.

The main idea is to show that even though $(M,\tilde{g})$ is compact, the $s$ parameter defined in \textbf{Lemma \ref{lemma_riccati}} goes to infinity, a phenomenon noted in \cite{BF}. Consequently, the $\tilde{g}$-distance to the boundary can still be viewed as a distance from infinity. Therefore, by taking a limit, we can derive super-harmonic functions. Indeed, adopting this viewpoint allows the construction of numerous super-harmonic functions.

The new Riccati equation can also be applied to static spaces with positive scalar curvature. We present the following theorem:
\begin{theorem}\label{thm_compact}
    Let $(M, g, V)$ be a static triple with positive scalar curvature. Suppose $\partial M = \Sigma \neq \emptyset$, $V = 0$ on $\Sigma$, and $V > 0$ elsewhere. Then the universal covering of $(M, g)$ is compact.
\end{theorem}
It is natural to consider the triple described in the theorem, as they can be viewed as foundational elements for constructing complete static triples. These manifolds exhibit numerous properties; for instance, $\Sigma$ is totally geodesic, and $|\nabla V|$ remains constant along $\Sigma$. A brief proof of these properties will be provided (also cf. \cite{Bo} and \cite{FM}). Consequently, by assembling these basic elements, we can construct larger static triples.

\textbf{Theorem \ref{thm_compact}} was first proved in \cite{Am} for 3-manifolds by lifting $(M, g)$ to a 4-manifold with constant Ricci curvature, followed by a Bonnet-Myers type argument. However, this approach can lead to singularities, necessitating careful handling. The new proof presented in this paper, however, does not encounter this issue.\\

The classification of 3-dimensional static triples is a current area of intense research. Let $g_0$ be the standard metric on $\mathbb{S}^2$. The only three known examples are:\\

     1) Standard round hemisphere: 
     
    \quad $(\mathbb{S}^3_+, g=dr^2+sin^2(r)g_0$ for $r\in [0,\frac{\pi}{2}], V=x_{n+1})$;
    
   2) Cylinder: 
    $([0,\frac{\pi}{\sqrt{3}}]\times\mathbb{S}^2,g=dt^2+g_{0},V=\sin(\sqrt{3}t)$;
    
    3) deSitter-Schwarzschild space: 
    
    \qquad $([r_1,r_2]\times\mathbb{S}^2,\frac{dr^2}{1-r^2-\frac{2m}{r}}+r^2g_{0},V=\sqrt{1-r^2-\frac{2m}{r}})$
    
    \qquad where  $m\in(0,\frac{1}{3\sqrt{3}})\text{ and }r_1<r_2 \text{ are positive zeros for } V$.\\
    
\begin{comment}
\begin{equation}
\begin{aligned}\label{example}
    &\text{1) Standard round hemisphere: }
    (\mathbb{S}^3_+, g=dr^2+sin^2(r)g_0 \text{ for } r\in [0,\frac{\pi}{2}], V=x_{n+1});\\
    &\text{2) Cylinder: }
    ([0,\frac{\pi}{\sqrt{3}}]\times\mathbb{S}^2,dt^2+g_{0},V=\sin(\sqrt{3}t);\\
    &\text{3) deSitter-Schwarzschild space:} \\
    &\qquad([r_1,r_2]\times\mathbb{S}^2,\frac{dr^2}{1-r^2-\frac{2m}{r}}+r^2g_{0},V=\sqrt{1-r^2-\frac{2m}{r}})\\
    &\qquad\text{where } m\in(0,\frac{1}{3\sqrt{3}})\text{ and }r_1<r_2 \text{ are positive zeros for } V.
    \end{aligned}
\end{equation}
\end{comment}
The locally flat static triples are well classified into these three categories, as detailed in \cite{Ko} and \cite{La}. We will show that the three static triples mentioned above correspond to solutions of the ordinary differential equation $\frac{\partial \theta}{\partial s} = -\frac{1}{n-1} \theta^2$ with different boundary conditions, potentially providing new insights into the classification problem.

Finally, we apply an integral identity to establish the following:

\begin{theorem}\label{thm-vanishing}
    Let $(M,V,g)$ be as in \textbf{Theorem \ref{thm_compact}}. Then $H^1(M) = 0$.
\end{theorem}

This paper is organized as follows. Section 2 provides the preliminaries. In Section 3, we establish \textbf{Lemma \ref{lemma_riccati}}, the central lemma of this study. Section 4 demonstrates the use of Busemann functions—distance functions parallel to the boundary—to prove \textbf{Theorem \ref{thm_splitting}}. Section 5 applies the central lemma to distance functions towards the boundary, leading to the proof of \textbf{Theorem \ref{thm_connect}} and a discussion of some related results. Finally, in Section 6, we extend these methodologies to positively curved static spaces.

\begin{comment}
\textbf{Acknowledgements}

The author is thankful to Prof Xiaodong Wang for bringing this problem to my attention, Prof Lucas Ambrozio, Demetre Kazaras, Greg Galloway, and Eric Woolgar for helpful discussions, and Prof Stefano Borghini, Mattia Fogagnolo and Chao Xia for mentioning 
\cite{BF} to me.

\textbf{Data Availability Statement}
The authors declare that no datasets were generated or analyzed during the current study. This manuscript does not include any associated data.

\textbf{Conflict of Interest Statement}
I hereby declare that I have no financial, personal, or professional conflicts of interest related to this work.
\end{comment}

\section{Preliminary}
In this section, we present basic definitions and related properties. The two propositions discussed in this section are well-known. Throughout this paper, terms with a tilde refer to those with respect to $\tilde{g} = \frac{1}{V^2}g$, and terms with an asterisk refer to those related to its universal covering.

\begin{proposition}
    Let $(M,g,V)$ be a static triple. Then it has constant scalar curvature.
\end{proposition}
\begin{proof}
    Taking $div_{g}$ to (\ref{static_equ_1}) gives
    \begin{equation}
            0=\frac{V}{2}dR+Ric(\nabla V,\cdot)+d(\Delta V)-div(\nabla^2 V)=\frac{V}{2}dR\notag
    \end{equation}
    We used contracted second Bianchi identity in the first equality and Ricci identity for the second equality.
\end{proof}
If we scale the metric so that the scalar curvature \( R  = \epsilon n(n-1) \), where \( \epsilon = -1, 0, 1 \), then taking the trace of (\ref{static_equ_1}) yields that (\ref{static_equ_1}) is equivalent to the following:
\begin{equation}
    \begin{aligned}
        V \left( \text{Ric} - n \epsilon g \right) - \nabla^2 V &= 0, \\
        \Delta V + n \epsilon V &= 0. \label{static_equ_2}
    \end{aligned}
\end{equation}

Set $\tilde{g}=\frac{1}{V^2}g$ and $\tilde{V}=\frac{1}{V}$. Under conformal change, the Ricci curvature term and Hessian terms are related by
\begin{equation}\notag
    \begin{aligned}
        Ric&=\tilde{Ric}-(n-2)\big(\frac{1}{V}\tilde{\nabla}^2 V-\frac{2}{V^2}dV\otimes dV\big) - \big(\frac{1}{V}\tilde{\Delta}V+\frac{n-3}{V^2}|\tilde{\nabla}V|_{\tilde{g}} \tilde{g}\big)\\
        \nabla^2V&=\tilde{\nabla}^2V+\frac{|\tilde{\nabla}V|_{\tilde{g}}}{V}\tilde{g}-\frac{2}{V}dV\otimes dV 
    \end{aligned}
\end{equation}
Plugging these into the static equation (\ref{static_equ_1}) or (\ref{sub_static_equ_1}), we obtain
\begin{equation}
   \tilde{V} \tilde{Ric}+(n-1)\tilde{\nabla}^2 {\tilde V}=(\geq)0\label{static_equ_tilde}
\end{equation}
respectively. We are going to use this inequality to derive Riccati type equation.\\

For negatively curved static spaces, we usually consider conformally compactifiable manifolds:
\begin{definition}
    A non-compact complete Riemannian manifold $(M,g)$ is called conformally compactifiable if $M$ is the interior of a compact manifold $\bar{M}$ with non-empty boundary $\Sigma$ and there exists a boundary defining function $r$ satisfying $r=0$ on $\Sigma$, $|dr|\neq 0$ on $\Sigma$, $r>0$ in $\mathring{M}$ so that $r^2g$ extends to a smooth metric on $\bar{M}$. $(\Sigma, r^2g\big|_{\Sigma})$ is called the conformal boundary.
\end{definition}
For static triples $(M,g,V)$, we usually assume that $\frac{1}{V}$ is the boundary defining function.\\

For positively curved static spaces, we have the following basic results.
\begin{proposition}\label{prop2}
    Let $(M,g,V)$ be as in \textbf{Theorem \ref{thm_compact}}. Then $|\nabla V|$ is a non-zero constant on $\Sigma$ and $\Sigma$ is totally geodesic.
\end{proposition}
\begin{proof}
    Since $\Delta V=-nV\leq 0$, Hopf lemma implies that $\frac{\partial V}{\partial \nu}<0$ where $\nu$ is outer normal vector. Note that $\Delta V=0$ on $\Sigma$ by the second equality in (\ref{static_equ_2}), therefore $\nabla^2V(X,Y)=\nabla^2V(X,\nu)=0$ on $\Sigma$ for any vectors $X,Y$ along $\Sigma$. And the proposition follows.
\end{proof}

\section{Proof of Riccati Equation}
In this section, we will prove \textbf{Lemma \ref{lemma_riccati}}. While the original proof is available in \cite{BF}, I will present it here due to its central importance to this study.
\begin{proof}
   Let \(\tilde{r}\) be a \(\tilde{g}\)-distance function, and let \(\tilde{A}\) and \(\tilde{H}\) denote the second fundamental form and mean curvature for the level sets of \(\tilde{r}\) with respect to \(\tilde{g}\). Under the conformal change, these terms, when computed with respect to \(g\), are given by
\begin{equation}
    \begin{aligned}\label{A_conformal}
        A&=V(\tilde{A}+\frac{\dot{V}}{V}\tilde{g})\\
        H&=\frac{1}{V}(\tilde{H}+(n-1)\frac{\dot{V}}{V})
    \end{aligned}
\end{equation}
where the dot  denotes differentiation along $\tilde{\nabla}\tilde{r}$. Set
\begin{equation}
    \theta\coloneqq \frac{H}{V}=\frac{1}{V^2}(\tilde{H}+(n-1)\frac{\dot{V}}{V})
\end{equation}
and take derivative of $\theta$, we get 
\begin{equation}
    \begin{aligned}\label{3_middle_1}
    \dot{\theta}&=-2\frac{\dot{V}}{V^3}(\tilde{H}+(n-1)\frac{\dot{V}}{V})+\frac{1}{V^2}\big[ \dot{\tilde{H}} +(n-1)\frac{\ddot{V}}{V}-(n-1)\frac{\dot{V}^2}{V^2} \big]\\
        &=-2\frac{\dot{V}}{V}\theta+\frac{1}{V^2}\big[ -\tilde{Ric}(\tilde{\nabla} \tilde{r},\tilde{\nabla} \tilde{r})-|\tilde{A}|_{\tilde{g}} +(n-1)\frac{\ddot{V}}{V}-(n-1)\frac{\dot{V}^2}{V^2} \big]
    \end{aligned}
\end{equation}
where we used Riccati equation for $\tilde{g}$ in the second equality. Next we want to use (sub)-static equation to get rid of $\tilde{Ric}$ term.
\begin{equation}
    \begin{aligned}
        \tilde{Ric}(\tilde{\nabla} \tilde{r},\tilde{\nabla} \tilde{r})&\geq-\frac{n-1}{\tilde{V}}\tilde{\nabla}^2\tilde{V}(\tilde{\nabla} \tilde{r},\tilde{\nabla} \tilde{r})\\
        &=-(n-1)V\ddot{(\frac{1}{V})}\\
        &=(n-1)\frac{\ddot{V}}{V}-2(n-1)\frac{\dot{V}^2}{V^2}\notag
    \end{aligned}
\end{equation}
Plugging this into (\ref{3_middle_1}), the \(\ddot{V}\) terms cancel out, and (\ref{3_middle_1}) becomes
\begin{equation}
    \dot{\theta}\leq -2\frac{\dot{V}}{V}\theta+\frac{1}{V^2}\big[(n-1)\frac{\dot{V}^2}{V^2}-|\tilde{A}|^2_{\tilde{g}}\big]\label{3_middle_2}
\end{equation}
From (\ref{A_conformal}) we know that $\tilde{A}=\frac{1}{V}A-\frac{\dot{V}}{V^3}g$, and it follows that
\begin{equation}
    |\tilde{A}|_{\tilde{g}}=V^2|A|_{g}^2-2\dot{V}H+(n-1)\frac{\dot{V}^2}{V^2}
\end{equation}
Plugging it into (\ref{3_middle_2}), we get
\begin{equation}
    \dot{\theta}\leq-|A|_g^2\label{3_middle_3}
\end{equation}
Finally, we want to make pick another parameter $s$ so that
\begin{equation}
    ds=V^2d\tilde{r}\label{s_parameter}
\end{equation}
In the new parameter, we get
\begin{equation}
    \begin{aligned}
        \frac{\partial \theta}{\partial s}&\leq -\frac{1}{V^2}|A|^2_g\\
        &\leq -\frac{1}{n-1}\frac{H^2}{V^2}=-\frac{1}{n-1}\theta^2\label{ODE}
    \end{aligned}
\end{equation}
We used Cauchy inequality in the second line.
\end{proof}

\begin{corollary}\label{cor_central}
    Under the same assumptions as in \textbf{Lemma \ref{lemma_riccati}}, let \(\gamma(t)\), \(t \in [0, T]\), be a \(\tilde{g}\)-minimizing geodesic with \(\gamma(0) = p\). Then, for any \(t \in (0, T)\),
    \begin{equation}
        \label{comparison}
        \theta(\gamma(t)) \leq \frac{n-1}{s(t)},
    \end{equation}
    where \(s(t) = \int_0^t V(\gamma(r))^2 \, dr\).
\end{corollary}

\begin{proof}
    Since \(\gamma\) is a minimizing geodesic in \([0, T]\), the distance function is smooth for \(t < T\). Thus, \(\theta\) is well-defined at \(\gamma(t)\).  For $0<t_2<T$, if $\theta(t_2) \leq 0$, then (\ref{comparison}) obviously holds at $t_2$. Hence, we may assume \(\theta(t_2) > 0\). From (\ref{ODE}), we know that \(\theta\) is non-increasing in $s$, thus in $t$, so \(\theta(t) > 0\) for $t\leq t_2$. From (\ref{ODE}) we know that $\frac{\partial}{\partial s}\frac{1}{\theta}\geq \frac{1}{n-1}$. Integrate it, we get
    \begin{equation}
        \frac{1}{\theta(t(s_2))}\geq \frac{s_2-s_1}{n-1}+\frac{1}{\theta(t(s_1))}\label{3_middle_5}
    \end{equation}
 Given that both sides of (\ref{3_middle_5}) are positive, it follows that
\begin{equation}
    \theta(t(s_2)) \leq \frac{\theta(t(s_1))}{1 + \left(\frac{s_2 - s_1}{n-1}\right) \theta(t(s_1))}.
\end{equation}
Thus, the corollary follows.

\end{proof}
(\ref{Riccati}) is equivalent to the Raychaudhuri equation in Lorentz metric. Lift $(M,g,V)$ and $(M,\tilde{g},V)$ to a $(n+1)$ spacetime $(N,h)$ and $(N,\tilde{h})$ respectively
\begin{equation}
    \begin{aligned}
        N&=\mathbb{R}\times M\\
        h&=-V^2dt^2+g, \quad \tilde{h}=\frac{1}{V^2}h=-dt^2+\tilde{g}\label{Lorentz_metric}
    \end{aligned}
\end{equation}
Given a unit-speed \(\tilde{g}\)-geodesic \(\gamma(t)\), we can lift it to a null geodesic \(\eta(t) = (t, \gamma(t))\) in \((N, \tilde{h})\), and \(t\) is the affine parameter, i.e., \(\nabla^{\tilde{h}}_X X = 0\) with \(X = \frac{\partial \eta}{\partial t}\). Under the conformal change, we can compute in \(h\) that
\begin{equation}
    \nabla^h_{X}X=\nabla^{\tilde{h}}_XX+2\frac{X(V)}{V}X-\frac{h(X,X)}{V}\tilde{\nabla}V=2\frac{X(V)}{V}X
\end{equation}
and it's still along $\frac{\partial \eta}{\partial t}$. This implies that \(\eta(t)\) remains a geodesic in \(h\), but \(t\) is no longer an affine parameter. If we define a new parameter \(s\) satisfying \(ds = V^2 dt\), and set \(\beta(s) = \eta(t(s))\), then direct computations show that \(\nabla^h_{\dot{\beta}} \dot{\beta} = 0\). This explains the choice of the new parameter in \textbf{Lemma \ref{lemma_riccati}}.

Let $\{e_i\}_{i=1}^{n}$ and $e_1=\frac{\dot{\gamma}}{V}$ be an ON frame along $\dot{\gamma}$ in $(M,g)$, with $\partial t$ they form a frame for $(N,h)$ along $\eta$. Let $K=\dot{\beta}$, and $\Pi$ be the projection $TN\big|_{p}\rightarrow (TN\big|_{p})/K$ for $p\in \beta$, i.e. $X,Y\in T_pN$ are equivalent iff $X-Y=cK$ for some $c$. Then the Weigarten map can be defined as
\begin{equation}
    \begin{aligned}
        b: (TN\big|_{p})/K&\rightarrow (TN\big|_{p})/K\\
        b([X])&=[\nabla^h_XK]
    \end{aligned}
\end{equation}
Direct computation shows that this map is well defined for vectors along $M$. Let $B=\sum_{i=2}^nh(b(e_i),e_i)$, then Raychaudhuri equation implies that
\begin{equation}
    B'\leq Ric_h(K,K)-\frac{1}{n-1}B^2\label{Raychaudhuri}
\end{equation}
If the triple $(M,g,V)$ is static, then $Ric_h(X,X)=0$ for all null vectors, which is called null energy condition. This inequality is equivalent to (\ref{Riccati}). For more details about Weigarten map and Raychaudhuri equations, readers are refered to \cite{EH}.

Now we have a better understanding for \textbf{Lemma \ref{lemma_riccati}}. We start with geodesics in $(M,\tilde{g})$, because they could be lifted to geodesics in $(\mathbb{R}\times M,-dt^2+\tilde{g})$ via $\eta(t)=(t,\gamma(t))$ where $\gamma(t)$ is a $\tilde{g}$ geodesic. It's important to note that geodesics in $(M,g)$ can't be lifted to geodesics in $(\mathbb{R}\times M,h=-V^2dt^2+g)$. We computes mean curvature for $g$ because $Ric_h(X,X)=0$ for null vectors so that we could apply Raychaudhuri equation. New parameter \(s(t) = \int_0^t V(\gamma(r))^2 \, dr\) is used because $s$ is the affine parameter for the curve $\eta$ in $h$.

\section{Splitting Theorem for Static Conformally Compactifiable Manifolds}
\begin{comment}
In this section we prove \textbf{theorem \ref{thm_splitting}}. Basically the proof is a mimic of the splitting theorem for complete Riemannian manifold of non-negative Ricci curvature in \cite{CG}. In the setting of non-negative Ricci curvature, we have the standard Ricci comparison result (\ref{Ricci_comparison}), and this shows that Busemann functions are super-harmoninc. Together with maximal principle, we get the desired splitting. In our case, \textbf{lemma \ref{comparison}} plays the role of (\ref{Ricci_comparison}).
\end{comment}
In this section, we prove \textbf{Theorem \ref{thm_splitting}}. The proof essentially mimics the splitting theorem for complete Riemannian manifolds with non-negative Ricci curvature as discussed in \cite{CG}. In scenarios of non-negative Ricci curvature, the standard Ricci comparison result (\ref{Ricci_comparison}) indicates that Busemann functions are super-harmonic. Utilizing the maximal principle, we achieve the desired splitting. In our case, \textbf{Lemma \ref{lemma_riccati}} serves the role of (\ref{Ricci_comparison}).

The outline of the proof is as follows: \textit{i)} We construct a minimizing geodesic line \(\gamma\) in \((M^*, \tilde{g}^*)\); \textit{ii)} then, we define a Busemann function associated to \(\gamma\) and demonstrate that they are super-harmonic, and thus harmonic by the maximum principle. \textit{iii)} we identify the desired splitting and iterate.

Throughout the proof, terms with \(*\) are those with respect to the universal covering, and terms with \(\tilde{}\) above are with respect to the conformal compactification.
\begin{proof}
\textit{i) Construct a minimizing line}
\begin{comment}
    Without loss of generality we might assume $(M^*,\tilde{g}^*)$ is not compact, otherwise there's nothing to prove. We could find points $p_i, q_i$ bounded away from the boundary so that $dis_{\tilde{g}}(p_i,q_i)\rightarrow \infty$. Let $\gamma_i$ be the minimizing geodesic connecting $p_i, q_i$ and $o_i$ be the middle points. Now fix a fundamentla domain $D$ for the universal covering $M^*$. By applying covering space transformations, we might assume all the $o_i$'s lie in $D$. Since $p_i$'s and $q_i$'s are universal bounded away from the boundary, \textbf{Condition C} implies that $\gamma_i$ are bounded away from the boundaries. By compactness of $\bar{D}$, we could find a subsequence so that $\gamma_i$ converges to a minimizing geodesic line $\gamma: \mathbb{R}\rightarrow M^*$ and $\gamma$ is universally away from boundary.\\
\end{comment}

    Without loss of generality, we may assume that \((M^*, \tilde{g}^*)\) is not compact; otherwise, there is nothing to prove. We can find points \(p_i, q_i\) bounded away from the boundary such that \(\text{dist}_{\tilde{g}}(p_i, q_i) \rightarrow \infty\). Let \(\gamma_i\) be the minimizing geodesics connecting \(p_i\) and \(q_i\), and let \(o_i\) be their midpoints. Now, fix a fundamental domain \(D\) for the universal covering \(M^*\). By applying covering space transformations, we may assume all the \(o_i\)'s lie in \(D\). Since \(p_i\)'s and \(q_i\)'s are uniformly bounded away from the boundary, \textbf{Condition C} implies that \(\gamma_i\) are bounded away from the boundaries. By the compactness of \(\overline{D}\), we can find a subsequence such that \(\gamma_i\) converges to a minimizing geodesic line \(\gamma: \mathbb{R} \rightarrow M^*\), and \(\gamma\) is universally away from the boundary.

\noindent \textit{ii) Busemann functions and superharmonicity}

Given a minimizing geodesic ray $c$, the associated Busemann function is defined as     
\begin{equation}
    \begin{aligned}
        b_c: M^*&\rightarrow \mathbb{R}\\
        b_c(x)&=\lim_{t\rightarrow \infty}(dis_{\tilde{g}}(x,\gamma_c(t))-t)\label{Busemann}
    \end{aligned}
    \end{equation}
   Using triangle inequality, it's easy to check that $dis_{\tilde{g}}(x,\gamma_c(t))-t$ is non-increasing in $t$ and bounded for fixed $x$, thus $b_{c}$ is well defined. Let $\gamma_+=\gamma\big|_{[0,\infty)}$ and $\gamma_{-}(t)=\gamma(-t)$ for $t\in [0,\infty)$, and $b_{\pm}$ be the Busemann functions associated to $\gamma_{\pm}$ respectively.

    Define the linear elliptic operator
\begin{equation}
\begin{aligned}\label{L_operator}
L\phi :=&\, \frac{1}{V^{*2}}\left [ \tilde{\Delta}^*\phi +\frac{(n-1)}{V^*}{\tilde g}^* ({\tilde \nabla}^*V^*,{\tilde \nabla}^*\phi)\right ]\\
=&\, {\tilde V}^{*2} \left [ \tilde{\Delta}^*\phi-\frac{(n-1)}{\tilde V^*}{\tilde g}^* ({\tilde \nabla}^*{\tilde V^*},{\tilde \nabla}^* \phi)\right ].
\end{aligned}
\end{equation}
This operator has the following geometric interpretation. Let \(\phi\) be a distance function in \(\tilde{g}\), i.e., \(|\tilde{\nabla} \phi|_{\tilde{g}} \equiv 1\). Then, the mean curvature of the level sets for \(\phi\) in \(\tilde{g}\) is \(\tilde{H} = \Delta_{\tilde{g}} \phi\). After a conformal change, its mean curvature in \(g\) is
\begin{equation}
    H = \tilde{V}^* \left(\tilde{\Delta}^*\phi-\frac{(n-1)}{\tilde V^*}{\tilde g}^* ({\tilde \nabla}^*{\tilde V^*},{\tilde \nabla}^*\phi)\right).
\end{equation}
Thus, \(V^* L(\phi) = H\) represents the mean curvature in \(g\) and \(L(\phi)\) is \(\theta\) as defined in \textbf{Lemma \ref{lemma_riccati}}. Next, we will show that \(L(b_{\pm}) \leq 0\) using \textbf{Corollary \ref{cor_central}}. Without loss of generality, let's work with \(b_+\). The remainder of this part follows closely the standard splitting proof, and only an outline will be provided (cf. \cite{Pe} for further details).

For any $p\in M^*$, we can construct a ray $\gamma_p$ starting from $p$ and ``parallel" to $\gamma_+$: for any $T\in [0,\infty)$, join $p$ with $\gamma_+(T)$ with unit velocity minimizing geodesic $\eta_T(t)$. Let $T\rightarrow \infty$, this sequence $\eta_T$ will subconverge to a minimizing ray, called $\gamma_p$, by the compactness of the set $\{X\in T_pM^*:\, |X|_{\tilde{g}}=1\}$. Such $\gamma_p$'s are called asymptotes for $\gamma_+$ from $p$. Let $b_p$ be the Busemann function associated to $\gamma_p$, then $b_p(x)+b_{+}(p)$ supports $b_{+}(x)$ from above at $p$, i.e.
\begin{lemma}
    \begin{equation}
        \begin{aligned}
            b_+(x)&\leq b_p(x)+b_+(p)\quad \forall x\in M^*\\
            b_+(p)&=b_p(p)+b_+(p)\notag
        \end{aligned}
    \end{equation}
\end{lemma}
See section 9.3 in \cite{Pe}, for example. As a result $L(b_+)(p)\leq L(b_p)(p)$ and it suffices to prove that $L(b_p)\leq 0$ at $p$. In order to show this, we find that $b_t(x)\coloneqq dis_{\tilde{g}}(x,\gamma_p(t))-t$ are support functions for $b_p(x)$ at $p$ for all $t>0$. Use \textbf{Corollary \ref{cor_central}}, we have that 
\begin{equation}
    L(b_p)(p)\leq L(b_t)(p)\leq \frac{n-1}{s(t)}\label{4_middle_1}
\end{equation}
where $s(t)=\int_0^tV(\gamma_p(r))^2dr$. Since \(\gamma_p\) is a minimizing ray, we allow \(t \to \infty\). Given that \(V\) is bounded from below, it follows that \(s(r) \to \infty\), and thus the right-hand side of (\ref{4_middle_1}) approaches zero. Consequently, it follows that \(L(b_+)(p) \leq 0\) for any \(p \in M^*\).

Now that we have established \(L(b_{\pm}) \leq 0\), it follows that \(L(b_+ + b_-) \leq 0\). Utilizing the triangle inequality, we find that \(b_+ + b_- \geq 0\) with equality along \(\gamma\). Since this superharmonic function achieves interior minimum, the maximum principle implies that \(b_+ + b_-\) is constant in \(M^*\). Therefore, we have \(0 \geq L(b_+) = -L(b_{-}) \geq 0\), which implies that \(L(b_{\pm}) = 0\).

\noindent \textit{iii)Splitting}\\
\begin{comment}
For simplicity, let $r=b_+$ in this section. Consider the level sets for $r$, and let $H$ be its mean curvature. We have shown that $L(r)$ can be interpreted as the $H/V$. Apply (\ref{Riccati}) to $\theta=H/V\equiv 0$, and equality holds in (\ref{Riccati}). Recall that we applied Cauchy inequality, so this implies that $A\equiv 0$ where $A$ is the second fundamental form w.r.t $g$. This implies that the level sets of $r$ remains the same. It follows that $g=u^2dr^2+h$ where $h$ is independent of $r$. It can be shown that $|\tilde{\nabla} r|_{\tilde{g}}=1$ since it's the limit of $\tilde{g}$ distance functioins, so $|\nabla b_+|_{g}=V$, so the $u=V^*$ in the splitting above, i,e,
\end{comment}
For simplicity, let \( r = b_+ \) in this section, and $r$ is a $\tilde{g}$ distance function. Consider the level sets for \( r \), and let \( H \) be its mean curvature. We have shown that \( L(r) \) can be interpreted as \( \frac{H}{V} \). Applying (\ref{Riccati}) to \(\theta = \frac{H}{V} \equiv 0\), we find that equality holds in (\ref{Riccati}). Recall that we applied the Cauchy inequality, so equality forces \( A \equiv 0 \) where \( A \) is the second fundamental form with respect to \( g \). This indicates that the level sets of \( r \) are totally geodesic. It follows that
    $g = u^2 dr^2 + h,$ where \( h \) is independent of \( r \). Since \( |\tilde{\nabla} r|_{\tilde{g}} = 1 \), being the limit of \(\tilde{g}\)-distance functions, we have \( |\nabla r|_{g} = V^*z \), thus \( u = V^* \) in the above splitting, i.e.
\begin{equation}
\begin{aligned}
M^*&=\mathbb{R}\times \Sigma\\
    g&=V^{*2}dr^2+h\label{splitting_middle}
    \end{aligned}
\end{equation}
    It remains to show that $V^*$ is independent of $r$. Let $\Sigma_r$ be level sets for $r$, and $\{\partial_i\}_{i=1}^{n-1}$ coming from local normal coordinates for $\Sigma$. Since $\Sigma_r$ are totally geodesic, by Codazzi equation we know that $Ric_{g^*}(\partial_i,\partial_r)=0$. Using static equation (\ref{static_equ_1}), we find that $\nabla^2_{g^*}V^*(\partial _r,\partial_i)=0$. It follows that
         \begin{equation}
     \begin{aligned}
         0=\nabla^2_{g^*}V^*(\partial_r,\partial_i)&=\partial_i\partial_rV^*-(\nabla_{\partial_i}\partial_r)V^*\\
         &=\partial_i\partial_rV^*-\frac{\partial_iV^*\partial_rV^*}{V^*}=V^*\partial_i\partial_r \log(V^*)\notag
     \end{aligned}
     \end{equation}
     This implies that $V^*$ splits as $V^*(r,y)=\alpha(r)\beta(y)$ for $y\in \Sigma$. Now examine $S$ in (\ref{static_equ_1}) along $\Sigma$ direction, i,e, $S(\partial_i,\partial_i)=0$. For $p\in M^*$, find local normal coordinates $\{x_i\}_{i=1}^{n-1}$ for $\Sigma$. Direct computation shows that
    \begin{equation}\label{connection}
        \nabla_{\partial r}\partial r=-V^*V^*_i\partial_i+\frac{V^*_r}{V^*}\partial r;\quad \nabla_{\partial_i}\partial_r=\nabla_{\partial_r}\partial_i=\frac{V^*_i}{V^*}\partial r;\quad \nabla_{\partial_i}\partial_j=0.
    \end{equation}
    Using these we can compute
    \begin{equation}
        \begin{aligned}
             R_{g^*}(\partial_i,\partial r,\partial r,\partial_i)&=\langle\nabla_{\partial r}\nabla_{\partial i}\partial_i,\partial r\rangle-\langle\nabla_{\partial i}\nabla_{\partial r}\partial_i,\partial r\rangle\\
             &=-\langle\nabla_{\partial i}(\frac{V^*_i}{V^*}\partial r),\partial r\rangle\\
             &=-V^*V^*_{ii}=-\alpha^2\beta\beta_{ii}
             \label{3_computation_curvature_1}
        \end{aligned}
    \end{equation}
    Using Gauss equation we have
    \begin{equation}
       Ric_{g^*}(\partial_i,\partial_i)=Ric_{h}(\partial_i,\partial_i)+R_{g^*}(\partial_i,\frac{\partial r}{V^*},\frac{\partial r}{V^*},\partial_i)=Ric_{h}(\partial_i,\partial_i)-\frac{\beta_{ii}}{\beta}\label{4_Ric_computation}
    \end{equation}
    $\nabla^2_{g^*}V^*$ can be computed as
    \begin{equation}
    \begin{aligned}\label{4_hessian}
        \nabla_{g^*}^2V^*(\partial_i,\partial_i)&=\alpha\beta_{ii}\\
        \nabla_{g^*}^2V^*(\frac{\partial_r}{V^*},\frac{\partial_r}{V^*})&=\frac{1}{V^{*2}}\big(\partial_r^2V^*-(\nabla^{g^*}_{\partial_r}\partial_r)V^*\big)\\
        &=\frac{\alpha_{rr}}{\alpha^2\beta}+\frac{\alpha|\nabla_{h}\beta|^2}{\beta}-\frac{\alpha_r^2}{\alpha^3\beta}\\        
        \Rightarrow \Delta_{g^*}V^*&=\alpha(\Delta_{h}\beta+\frac{|\nabla_{h}\beta|^2}{\beta})+\frac{1}{\beta}(\frac{\alpha_{rr}}{\alpha^2}-\frac{\alpha_r^2}{\alpha^3})
        \end{aligned}
    \end{equation}
    Put these together, we get
    \begin{equation}
    \begin{aligned}
        0&=\sum_{i=1}^{n-1}\big(V^*Ric_{g^*}(\partial_i,\partial_i)+\Delta_{g^*}V^*-\nabla^2_{g^*}V^*(\partial_i,\partial_i)\big)\\
        &=\alpha\big[\beta R_{h}+(n-3)\Delta_{h}\beta+(n-1)\frac{|\nabla_{h}\beta|^2}{\beta})\big]+(n-1)\frac{1}{\beta}\big[\frac{\alpha_{rr}}{\alpha^2}-\frac{\alpha_r^2}{\alpha^3}\big]\\
       \Rightarrow & \beta\big[\beta R_{h}+(n-3)\Delta_{h}\beta+(n-1)\frac{|\nabla_{h}\beta|^2}{\beta})\big]=-(n-1)\frac{1}{\alpha}\big[\frac{\alpha_{rr}}{\alpha^2}-\frac{\alpha_r^2}{\alpha^3}\big]
    \end{aligned}
    \end{equation}
    This equation is seperable, which implies that $\frac{1}{\alpha}\big[\frac{\alpha_{rr}}{\alpha^2}-\frac{\alpha_r^2}{\alpha^3}\big]=C$. If we make the coordinate transformation $dt=\alpha dr$, this implies that
    \begin{equation}
        \alpha \frac{\partial^2 \alpha}{\partial t^2}=c \label{3_middle_4}
    \end{equation}
    Note that $V^*$ comes the lifting of $V$, which is bounded both below and above when away from the boundary. So $\alpha$ is also bounded. The only bounded solution to (\ref{3_middle_4}) are constant solutions, so $V^*$ doesn't depend on $r$, and we are done.\\

    \textit{iv) Iteration}
    
    Now, having proven the splitting (\ref{splitting_middle}), we observe that \( V^* \) does not depend on \( r \), allowing \( V^* \) to be viewed as a function on \( \Sigma \). If \((\Sigma, \frac{1}{V^{*2}}h)\) is not compact, as previously, we can construct a geodesic line \( \gamma_2 \). However, note that \((\Sigma, h, V)\) no longer forms a static triple, i.e. (\ref{static_equ_1}) doesn't hold, thus the prior argument does not facilitate a splitting for \((\Sigma, \frac{1}{V^{*2}}h)\). Instead, we consider \( \gamma' \) as a geodesic line in \((M, \tilde{g}^*, V^*)\) by identifying \( \Sigma \) as the slice \(\{0\} \times \Sigma\). Repeating the previous steps results in another splitting:
\begin{equation}\notag
    \tilde{g}^* = ds^2 + h_2.
\end{equation}
This new splitting differs from the one derived from \( \gamma \) because \( \partial_s \) and \( \partial_r \) are orthogonal to each other along \( \gamma_2 \subset \{0\} \times \Sigma \). By repeating this process, the analysis is completed.

    \end{proof}

If the triple \((M, g, V)\) is only sub-static, it is straightforward to see that \textit{i)} and \textit{ii)} of the proof still apply, and we achieve the splitting $g = V^{*2} \, dr^2 + h$.
Furthermore, we claim that \( V^* = \alpha(r) \beta(y) \). To establish this, we need to verify that \( S(\partial_r, \partial_r) = S(\partial_i, \partial_r) = 0 \). In fact, \( S(\partial_r, \partial_r) = 0 \) is a direct consequence of the splitting \( g = V^{*2} \, dr^2 + h \).
 Using first half of (\ref{3_computation_curvature_1}), we have that
    \begin{equation}
        \begin{aligned}
            S(\frac{\partial r}{V^*},\frac{\partial r}{V^*})=-\sum_{i=1}^{n-1}V^*_{ii}+\Delta V^*-\nabla^2V^*(\frac{\partial r}{V^*},\frac{\partial r}{V^*})=0
        \end{aligned}
    \end{equation}
    $V_{ii}=\nabla^2V(\partial_i,\partial_i)$ because we picked normal coordinate for $\Sigma$ and $\Sigma$ is totally geodesic in $M$. Next consider the function 
    \begin{equation}
        f(t)\coloneqq S(\frac{\partial_r}{V}+t\partial_i,\frac{\partial_r}{V}+t\partial_i)
    \end{equation}
     Since we have shown that $S(\frac{\partial r}{V},\frac{\partial r}{V})=0$, $f(t)=2S(\frac{\partial r}{V},\partial_i)t+S(\partial_i,\partial_i)t^2$, and $f(t)\geq 0$ for all $t\in \mathbb{R}$ by sub-static condition. This implies $S(\frac{\partial r}{V},\partial_i)=0$ by the determinant of the quadratic. Then the same argument gives the splitting $V^*=\alpha(r)\beta(y)$. However, it remains uncertain how \( S\big|_{\Sigma} \geq 0 \) influences the function \( \alpha \). Specifically, proving that \( \alpha \) is constant requires further exploration of the implications of \( S\big|_{\Sigma} \geq 0 \) on the curvature and other geometric properties of \( \Sigma \).\\

An important inquiry pertains to the possibility of omitting \textbf{Condition C} from \textbf{Theorem \ref{thm_splitting}}. In the original proof, we utilized this condition to construct $\tilde{g}$-geodesics distanced from the boundary. As discussed in the subsequent section, even if the geodesic line $\gamma$ interacts with the boundary, the associated Busemann functions $b_{\pm}$ remain super-harmonic. The critical failure occurs with the maximum principle. If any segment of $\gamma$ remains within $M^*$, the maximum principle can be invoked to achieve the desired splitting. However, there exists the possibility that $\gamma$ entirely resides within $\Sigma^*$.
\\

     If the conformal boundary is a flat torus, the following uniqueness result serves as a corollary in \cite{GSW}:
     \begin{corollary}
         Let $(M,g,V)$ be as \textbf{Theorem \ref{thm_splitting}}. Further assume that
         \begin{enumerate}
             \item The conformal boundary $(\Sigma,h)$ is flat $(n-1) $torus $\big(\mathbb{T}^{n-2}\times \mathbb{S}^1,\sum_{i=1}^{n-2}d\theta^2+d\phi^2)$, and the period of $\phi$ factor is $\frac{4\pi}{nr_0}$;
             \item Let $i$ be the inclusion map $\Sigma\rightarrow \tilde{M}$, then the kernel of the induced map between fundamantal groups $i_*:\Pi(\Sigma)\rightarrow \Pi(\tilde{M})$ is generated by $S^1$ factor;
             \item $\frac{\partial R_{\tilde{g}}}{\partial \nu^{n-2}}<0$ pointwise on $\Sigma$, where $\nu$ is inner normal derivative for $\tilde{g}$            
         \end{enumerate}
         Then $(M,g,V)$ is isometric to the Horowitz-Myers soliton given by
         \begin{equation}
         \begin{aligned}
         M&=\mathbb{R}^2\times \mathbb{T}^{n-2}\\
             g&=\frac{1}{r^2(1-\frac{r_0^n}{r^n})}dr^2+r^2(1-\frac{r_0^n}{r^n})d\phi^2+r^2\sum_{i=1}^{n-2}d\theta^2\\
             V&=r
             \end{aligned}
         \end{equation}
     \end{corollary}

\begin{remark}One reason why Horowitz-Myers soliton is important is that it's conjectured to be the ground state for Wang's mass, see \cite{Wo}. And the conjecture has recently been proved for $n\leq 7$ in \cite{BH}.

\end{remark}

\section{Connectedness of Conformal Boundary}

In this section, we will explore how \textbf{Lemma \ref{lemma_riccati}} can assist in understanding the topology of conformally compactifiable static manifolds, specifically through the super-harmonicity of distance functions and properties of minimal surfaces.

In the previous section, we constructed a \(\tilde{g}^*\) minimizing line that remains distanced from the conformal boundary and associated with it Busemann functions \(b_{\pm}\). These functions can be considered as distance functions from infinity. Using (\ref{lemma_riccati}) and taking a limit, we demonstrated that \(b_{\pm}\) are super-harmonic. In this section, we will instead use distance functions to the boundary. Interestingly, even though \(\tilde{g}\) is compact, due to the parameter transformation described in (\ref{parameter}), these distance functions to the boundary can still be regarded as distances from infinity. This property has been extensively analyzed in \cite{BF}.

     \begin{lemma}\label{riccati_boundary}
         Suppose \((M, g, V)\) is a complete conformally compactifiable sub-static triple with conformal boundary \(\Sigma\), and \(\tilde{V}\) is a boundary defining function. For any \(p \in \Sigma\), let \(u(x) = \text{dist}_{\tilde{g}}(x, p)\). Then, the operator \(L\), defined as in (\ref{L_operator}) but with \(g^*\) and \(V^*\) replaced by \(g\) and \(V\), satisfies \(L(u) \leq 0\) everywhere in the support sense.

     \end{lemma}
     \begin{proof}
         For any $q\in M$, connect $q$ with $p$ by a $\tilde{g}$-minimizing unit speed curve $\gamma(t)$ with $\gamma(0)=p$, $\gamma(T)=q$ where $T=dis_{\tilde{g}}(p,q)$. Note that $(\bar{M},\tilde{g})$ is a compact manifold with boundary $\Sigma$, and $\gamma$ might fail to be a geodesic when touching $\Sigma$. First assume $\gamma$ avoids all the other boundary component, then $\gamma$ is a minimizing geodesic. Let $v_{\epsilon}(x)=d_{\tilde{g}}(\gamma_1(\epsilon),x)$ for $\epsilon>0$, then $v_{\epsilon}(x)+\epsilon$ is a support function for $u_1(x)$ at $q$ from above, i.e.
    \begin{equation}
        \begin{aligned}\notag
            u_1(x)\leq v_{\epsilon}(x)+\epsilon\\
            u_1(q)= v_{\epsilon}(q)+\epsilon
        \end{aligned}
    \end{equation}
    Since $\gamma$ is a minimizing geodesic, $v_{\epsilon}(x)$ is smooth at $x$. And we have $L(u)(q)\leq L(v_{\epsilon})(q)$, so it suffices to show that $L(v_{\epsilon})(q)\rightarrow 0$ as $\epsilon\rightarrow 0$.
    
    By applying \textbf{Corollary \ref{cor_central}}, we have that
    \begin{equation}
        L(v_\epsilon)\leq \frac{n-1}{s(\epsilon)},\text{ with } s(\epsilon)=\int_{\epsilon}^{T}\frac{1}{\tilde{V}(\gamma(r))^2}dr
    \end{equation}
    It suffices to show that $\int_{\epsilon}^{T}\frac{1}{\tilde{V}(\gamma(r))^2}dr\rightarrow \infty$ as $\epsilon\rightarrow \infty$. Since \(\tilde{V}\) is a boundary defining function, \(|\tilde{\nabla} \tilde{V}|\) is non-zero on \(\Sigma\). Given that \(\tilde{V}\) extends smoothly to \(0\) on \(\Sigma\) and \((\bar{M},\tilde{g})\) is compact, we infer that $\tilde{V}(\cdot)\leq cd_{\tilde{g}}(\cdot,\Sigma_1)\leq cd_{\tilde{g}}(\cdot, p)$ for some constant $c>0$. Thus we have along $\gamma_1$ that $V(\gamma(r))\geq \frac{1}{cr}$.
\begin{equation}
    \int_{\epsilon}^{d_{\tilde{g}}(x,p)}V(\gamma(r))^2dr\geq \int_{\epsilon}^{d_{\tilde{g}}(x,p)}\frac{1}{c^2r^2}dr\rightarrow \infty
\end{equation}
This case is proved. Now suppose $\gamma$ touches boundary for at some $\gamma(t), 0<t<T$. Let $T_1=\sup\limits_{0<t<T}\{\gamma(t)\in\Sigma\}$, i.e. the last time $\gamma$ touches boundary. Then we could define $w(x)=dis_{\tilde{g}}(x,\gamma(T_1))$, and $\eta(r)=\gamma(T_1+r)$. Then we could apply the above argument to $w$ and $\eta$. $w(x)+T_1$ supports $u(x)$ from above, and we are done.
     \end{proof}

     \noindent \textit{Proof of theorem \ref{thm_connect}}
     
Suppose \(\partial \bar{M} = \Sigma\) has more than one connected component, with \(\Sigma_1\) and \(\Sigma_2\) as two of them. Fix points \(p_1 \in \Sigma_1\) and \(p_2 \in \Sigma_2\), and connect them with a length minimizing curve \(\gamma\) (might not be a geodesic on boundary). Let \(u_1(x) = \text{dist}_{\tilde{g}}(x, p_1)\) and \(u_2(x) = \text{dist}_{\tilde{g}}(x, p_2)\). By \textbf{Lemma \ref{riccati_boundary}}, we have that \(L(u_i) \leq 0\) for \(i = 1, 2\). Consequently, \(u_1 + u_2\) is a super-harmonic function and attains an interior minimum along \(\gamma\) due to the triangle inequality. By the maximum principle, \(u_1 + u_2\) must be constant. However, this leads to a contradiction since equality should only hold along \(\gamma\) by the triangle inequality, concluding the proof.

         \qed\\

\textbf{Lemma \ref{riccati_boundary}} provides a wealth of super-harmonic functions. We can also demonstrate that \(\text{dist}_{\tilde{g}}(x, \Sigma)\) is super-harmonic. Due to the blow-up of the \(s\) parameter as described in (\ref{parameter}), we can still refer to these functions as distances from infinity. However, for positively curved static triples as discussed in \textbf{Theorem \ref{thm_compact}}, despite \(\tilde{g}\) not being compact, the new parameter transformation described in (\ref{parameter}) remains finite.
\\
\begin{remark}
In \cite{BF}, the authors introduced the notion of \(V\)-completeness and constructed Busemann functions for \(V\)-complete ends. In our setting of conformally compactifiable manifolds, we can more concretely construct these Busemann functions as the \(\tilde{g}\) distance to points on the boundary, which simplifies the analysis.
\end{remark}

\textbf{Theorem \ref{thm_connect}} can be proved by applying a different interpretation of (\ref{Riccati}). Suppose that \(\Sigma\) has more than one connected component, and consider \(\Sigma_1\) as one of them. In this context, \(\Sigma_1\) constitutes a non-trivial element in \(H_{n-1}(M)\), the \((n-1)\)-th homology group of \(M\). According to geometric measure theory \cite{FF}, we can identify an interior stable minimal surface \(\Sigma'\) that represents the homology class \([\Sigma_1] \in H_{n-1}(M)\), i.e.
\begin{equation}
    H=0 \text{ and } \frac{\partial H}{\partial \nu}\geq 0\label{minimal_surface}
\end{equation}
where $\nu$ is unit vector normal to $\Sigma'$. However, by running the flow \(\Phi(t, x) = V \nu\), which corresponds to the geodesic flow in \(\tilde{g}\), we derive the inequality: $\frac{\partial \theta}{\partial s} \leq -\frac{\theta^2}{n-1}$ with equality if and only if \(A = 0\). Here, \(\theta\) and \(s\) retain their prior definitions. Alternatively, this relationship can also be observed by applying \(V\) within the stability operator framework (cf. \cite{Li}).

\begin{equation}
    \begin{aligned}
        0&\geq \int_{\Sigma} V\big(\Delta_{\Sigma'}V+V|A|^2+VRic_{g}(\nu,\nu)\big)\\
        &=\int_{\Sigma} V\big[\big(\Delta_{g}V-\nabla^2_gV(\nu,\nu)-Hg(\nabla V,\nu)\big)+V|A|^2+VRic_{g}(\nu,\nu)\big]\\
        &\geq\int_{\Sigma}V^2|A|^2\notag
    \end{aligned}
\end{equation}
where we used sub-static inequality (\ref{sub_static_equ_1}) and $H=0$ in the last inequality. As a result, the metric splits locally as $g=(V\circ \phi)^2dt^2+g\big|_{\Sigma'}$. Continuing to move along the flow \(\Phi\), it becomes evident that this splitting holds universally. However, such metrics are not conformally compactifiable. For a more detailed discussion on the computation, please refer to the argument presented in Proposition 19 in \cite{Am}.\\

\begin{comment}
If we assume $Ric\geq -(n-1)g$, we have \textbf{theorem \ref{Ricci_connectedness}} that shows the conformal boundary is connected if one boundary component has non-negative Yamabe invariant. The condition which requires non-negative Yamabe invariant can not be removed. In \cite{CG2}, M.Cai and G.Galloway used Riccati equation, but in $Ric\geq -(n-1)g$ case, the equation will give us
\begin{equation}
    H'\leq (n-1)-\frac{H^2}{n-1}\label{Riccai_Ric}
\end{equation}
and $H\equiv -(n-1)$ is a solution. If one boundary component has negative scalar curvature, we could construct a surface near boundary such that $H<1-n$. then we could apply (\ref{Riccai_Ric}) to get desired result. If Yamabe invariant equals $0$, we will have to involve Busemann function, and we'll not discuss it here. In our case, instead, the Riccati equation gives us (\ref{Riccati}) and we used it to construct super-harmonic functions, and that's why we are not assuming the boundary has components of non-negative Yamabe invariant.
\end{comment}
If we assume \( \text{Ric} \geq -(n-1)g \), \textbf{Theorem \ref{Ricci_connectedness}} demonstrates that the conformal boundary is connected if one boundary component has a non-negative Yamabe invariant. However, this condition requiring a non-negative Yamabe invariant cannot be removed. In \cite{CG2}, M. Cai and G. Galloway employed the Riccati equation. In the case of \( \text{Ric} \geq -(n-1)g \), the equation yields
\begin{equation}
    H' \leq (n-1) - \frac{H^2}{n-1}, \label{Riccai_Ric}
\end{equation}
with \( H \equiv -(n-1) \) being a particular solution. If a boundary component has negative scalar curvature, one can construct a surface near the boundary such that \( H < 1 - n \). Applying (\ref{Riccai_Ric}) would then lead to the desired result. If the Yamabe invariant is zero, the analysis would involve Busemann functions, which we do not cover here. In our approach, the Riccati equation provides (\ref{Riccati}) and is used to construct super-harmonic functions, thus circumventing the need for assuming boundary components with a non-negative Yamabe invariant.

Actually, we have the following example. Let $(\Sigma^{n-1},h)$ be a compact Riemannian manifold satisfying $Ric_{\Sigma}=-(n-2)h$. Then
\begin{equation}
(\mathbb{R}\times \Sigma, g=dt^2+\cosh^2(t)h)\label{exmaple_unconnected}
\end{equation}
is complete manifold with compactification $[0,1]\times \Sigma$ and $Ric_{g}=-(n-1)g$.

X. Wang \cite{Wa2} correlates \(\lambda_1\), the infimum of the \(L^2\) spectrum of \(-\Delta\), with the \((n-1)\)-th homology group by asserting that \(H_{n-1}(M, \mathbb{Z}) = 0\) provided \(\lambda_0 > n-2\). If \(\lambda_0 = n-2\) and \(H_{n-1}(M, \mathbb{Z}) \neq 0\), then \(M, g\) must correspond to the example discussed above.

\section{Applications to Static Triples of Positive Scalar Curvature}
In this section, we will explore the application of the Riccati equation to positively curved static triples and identify the three known examples as distinct initial value solutions to the ordinary differential equation $\theta' = -\frac{\theta^2}{n-1}$.\\

\noindent \textit{Proof of theorem \ref{thm_compact}}\\
\begin{comment}
It suffices to show that any $\alpha\in\Pi_1(M)$ must be of finite order, and we prove by contradiction. Suppose $\alpha$ is such an element, and let $\Gamma\in \Pi_1(M)$ generated by $\alpha$. Let $M^*$ be the Riemannian covering space w.r.t $\Gamma$.
\end{comment}
Let $\Sigma_{\epsilon}=\{x\in M:V(x)=\epsilon\}$. The second fundamental form $\tilde{A}$ and $A$ are related by\\
\begin{equation}
    \notag \tilde{A}=\frac{1}{V}(A+\frac{\partial }{\partial n}(-\log V)g\big|_{\Sigma_{\epsilon}})=\frac{1}{V}(A-\frac{\partial V}{\partial n}g\big|_{\Sigma_{\epsilon}})
\end{equation}
where $n$ is outer normal. By Proposition \ref{prop2}, the boundary $\partial M = \Sigma$ is totally geodesic, and $\frac{\partial V}{\partial n}$ is a negative constant on $\Sigma$. Consequently, $\tilde{A}$ is positive definite for $\Sigma_{\epsilon}$ when $\epsilon$ is sufficiently small. This corresponds to \textbf{Condition C} in \textbf{Theorem \ref{thm_splitting}}. Repeating the proof in that theorem, we find that, due to convexity, all geodesic lines and rays constructed as in \textbf{Theorem \ref{thm_splitting}} remain away from the boundary. Therefore, $s(t) = \int_0^t V^{*2} \, dr \rightarrow \infty$ along these lines and rays. The same argument applies to obtain the splitting result:
\begin{equation}
M^*=\mathbb{R}^k\times \Sigma^{n-k}\quad g^*=V^{*2}g_{\mathbb{E}}+h\label{middle_thm5}
\end{equation}
where $(\Sigma,h)$ is compact with non-empty boundary, and $k\geq 1$. Fix a point \( p \in \partial M \) and let \( U \) be a small neighborhood in $M$ surrounding \( p \). According to (\ref{middle_thm5}), the local metric within \( U \) splits as \( g = V^2 g_{\text{flat}} + h \). Since \( V = 0 \) on \( \partial M \), the term \( V^2 g_{\text{flat}} \) collapses, which implies that \( \partial M \) cannot locally behave as an \((n-1)\)-dimensional manifold, thereby leading to a contradiction.

\qed
\\
This result was also proved in \cite{Am} in dimension 3. The idea is to lift $(M,g,V)$ to a 4-dimensional manifold via
\begin{equation}
    (\mathbb{S}^1\times M, h=V^2d\phi^2+g)
\end{equation}
and $h$ has constant positive Ricci curvature, but possibly with singularities. Using classical Bonnet-Myers type argument and generalized Hopf-Rinow theorem for singular space, and the result follows. In this new proof, we avoid involving singularities. Moreover, the proof presented in \cite{Am} is specific to three-dimensional spaces because it depends on the compactness of the universal covering of the boundary. This requirement is supported by a result from \cite{SY}, which demonstrates that each stable component of \(\partial M\) is diffeomorphic to a sphere. This conclusion is feasible because \(M\) is oriented and exhibits positive scalar curvature in three-dimensional contexts.
\\

\textbf{Theorem \ref{thm_compact}} can be viewed as the positive curvature analog of \textbf{Theorem \ref{thm_splitting}}, focusing on geodesic lines parallel to the boundary. Moreover, there exists a corresponding result similar to \textbf{Theorem \ref{thm_connect}}, which employs the distance function to the boundary.

\begin{theorem}[L.Ambrozio]
    Let $(M,V,g)$ be as in \textbf{Theorem \ref{thm_compact}}. Then $\partial M=\Sigma$ can have at most one unstable boundary component.
\end{theorem}
In \cite{Am}, the proof is based on the fact that there doesn't exist interior minimal surfaces except cylinder, and an application of geometric measure theory \cite{FF}, as in the discussion in the alternative proof for \textbf{Theorem \ref{thm_connect}} presented in the previous section. If we strengthen the condition unstable minimal boundary components to that $H< 0$ for the level sets $\{V=\epsilon\}$ when $\epsilon$ is small, we could also apply the method in \textbf{Theorem \ref{thm_connect}} to prove this result. Since the new assumption is stronger than that in the theorem above and the proof is similar to that in \textbf{Theorem \ref{thm_connect}}, only a sketch will be given:\\

\noindent\textit{Sketch}:
Let \( H \) denote the mean curvature with respect to the inner normal. Suppose \( H < 0 \) for the level sets \(\{V = \epsilon\}\) for small \(\epsilon\) near boundary components \(\Sigma_1\) and \(\Sigma_2\). Let \(\Sigma_1'\) and \(\Sigma_2'\) be two such surfaces near \(\Sigma_1\) and \(\Sigma_2\), respectively. Define \(M'\) as the subset of \(M\) contained between \(\Sigma_1'\) and \(\Sigma_2'\). Let \(u_i(x) = \text{dist}_{\tilde{g}}(x, \Sigma_i')\) for \(i = 1, 2\), and let \(L(u_i)\) be defined as before. The condition \(H < 0\) implies that \(L(u_i) < 0\) on \(\Sigma_i\). By \textbf{Lemma \ref{lemma_riccati}}, this implies \(L(u_i) < 0\) in \(M'\) in the support sense. Consequently, \(u_1 + u_2\) is a strictly superharmonic function. However, \(u_1 + u_2\) achieves its infimum along the geodesic minimizing connecting \(\Sigma_1\) and \(\Sigma_2\), leading to a contradiction.\qed
\\

Another intriguing application of the Riccati equation is the classification of the three known examples in introduction. We will explain the classification process. Let \( g_0 \) denote the standard metric on \(\mathbb{S}^2\).

For hemisphere, we can write the metric as $g_{can}=dr^2+\sin(r)^2g_0$, $r\in[0,\pi/2]$, $V=\cos(r)$, and $\tilde{g}=\frac{dr^2}{\cos(r)^2}+\tan(r)^2g_0$.  Let $\Sigma_r$ be the level sets for $r$. The $g$-mean curvature for $\Sigma_{r}$ is $H=2\cot(r)$, and therefore $\theta=H/V=\frac{2}{\sin(r)}$. The $s$-parameter is given by $ds=\cos(r)^2dt=\cos(r)dr$ where $t$ is $\tilde{g}$ distance fucntion and thus we can pick $s=\sin(r)$. So they solve the following ODE
\begin{equation}
    \begin{aligned}\label{ODE_hemisphere}
        \frac{\partial \theta}{\partial s}&=-\frac{\theta^2}{2}\\
        \theta&=2 \text{ for } s=1,\text{ which is on }\Sigma\\
        \theta&\rightarrow \infty \text{ as }s\rightarrow 0, \text{ which is north pole }
    \end{aligned}
\end{equation}

For cylinder $([0,\frac{\pi}{\sqrt{3}}]\times\mathbb{S}^2,dr^2+g_0,V=\sin(\sqrt{3}r)$, apparently $\{r\}\times \mathbb{S}^{2}$ gives, with a reparametrization, $\tilde{g}$ geodesic flow, and $H=0$. Thus $\theta\equiv0$ sovles the ODE.

Finally for deSitter-Schwarzschild metric
\begin{equation}
     \Big([r_1,r_2]\times\mathbb{S}^2,\frac{dr^2}{1-r^2-\frac{2m}{r}}+r^2g_0,V=\sqrt{1-r^2-\frac{2m}{r}}\Big)
\end{equation}
where $m\in(0,\frac{1}{3\sqrt{3}})$ and $r_1<r_2$ are positive zeros for $V$. With a similar argument, we have $H=2\frac{\sqrt{1-r^2-\frac{2m}{r}}}{r}$ and $\theta=\frac{2}{r}$. The relationship between $r$, $g$-distance $l$, $\tilde{g}$-distance $t$ and $s$ parameter are
\begin{equation}
    ds=V^2dt=Vdl=dr\notag
\end{equation}
So we can take $s=r$. As a result, they solve the following
\begin{equation}
\begin{aligned}
    \frac{\partial \theta}{\partial s}&=-\frac{\theta^2}{2}\\
    \theta(r_1)&=\frac{2}{r_1}; \theta(r_2)=\frac{2}{r_2};
\end{aligned}
\end{equation}

From the discussion above, it is apparent that the three known positively curved static triples correspond to solutions of the ODE \(\frac{\partial \theta}{\partial s} = -\frac{\theta^2}{2}\), each with distinct boundary values. Additionally, the traceless part of the second fundamental forms associated with \(\theta\) vanishes in all cases. This observation may facilitate a complete classification of such geometries.
\\

Next, we construct a family of solutions to the static equation (\ref{static_equ_1}). From \textbf{Theorem \ref{thm_splitting}} and \textbf{Theorem \ref{thm_compact}}, we observe that the metric of the form
\begin{equation}
    M = S^1 \times \Sigma, \quad g = V^2 \, dr^2 + h \label{6_splitting}
\end{equation}
where \( V \) and \( h \) are independent of \( r \), is of particular interest. This is because (\ref{static_equ_1}) is automatically satisfied in the \((\partial_r, \partial_r)\) direction, as stated in the argument after the proof of \textbf{Theorem \ref{thm_splitting}}, and in the \((\partial_r, X)\) direction for \( X \) along \(\Sigma\) by the Codazzi equation. Consequently, (\ref{static_equ_1}) reduces to an equation solely on \((\Sigma, V, h)\). Identifying \(\Sigma\) as a slice \(\{p\} \times \Sigma\) within \(M\), \(V\) can be viewed as a function in both \(M\) and \(\Sigma\). Since \(\nabla_g^2 V \big|_{\Sigma} = \nabla_h^2 V\), and considering computations that demonstrate (\ref{4_Ric_computation}) and (\ref{4_hessian}), the static equation for the metric of the form (\ref{6_splitting}) simplifies to
\begin{equation}
    S' \coloneqq V \text{Ric}_h + \left(\Delta_h V + \frac{1}{V} |\nabla_h V|^2 \right) h - 2 \nabla_h^2 V = 0 \label{static_n-1}
\end{equation}

Now we further assume $(\Sigma,V,h)$ is of the form
\begin{equation}
    \begin{aligned}
        &\Sigma=I\times \mathbb{S}^{n-2};\quad V=V(s)\\
        &h=ds^2+f(s)^2g_0
    \end{aligned}
\end{equation}
where $I$ is an interval of $\mathbb{R}$ and $g_0$ is the round metric on $\mathbb{S}^{n-2}$. Let $X,Y$ be unit vectors along $\mathbb{S}^{n-2}$. We computes $Ric_{h}$ for the warped product:
\begin{equation}
    \begin{aligned}\label{6_Ric}
        Ric_h(\partial s,\partial s)&=-(n-2)\frac{\ddot{f}}{f}\\
        Ric_h(\partial_s,X)&=0\\
        Ric_h(X,Y)&=\big((n-3)\frac{1-\dot{f}^2}{f^2}-\frac{\ddot{f}}{f}\big)h(X,Y)\\
    \end{aligned}
\end{equation}
and $\nabla_h^2V$
\begin{equation}
    \begin{aligned}
        \nabla_h^2V(\partial s,\partial s)&=\ddot{V}\\
        \nabla_h^2V(\partial s,X)&=0\\
        \nabla_h^2V(X,X)&=\frac{\dot{f}\dot{V}}{f}\label{6_hessian}
    \end{aligned}
\end{equation}
Note that $S'(\partial_s,X)=0$ automatically, therefore we only need $S'(\partial_s,\partial_s)=0$ and $S'(X,X)=0$, which are
\begin{equation}
    \begin{aligned}
      & -(n-2)\frac{\ddot{f}V}{f}+\ddot{V}+(n-2)\frac{\dot{f}\dot{V}}{f}+\frac{\dot{V}^2}{V}-2\ddot{V}=0\\
       & V\big((n-3)\frac{1-\dot{f}^2}{f^2}-\frac{\ddot{f}}{f}\big)+\ddot{V}+(n-2)\frac{\dot{f}\dot{V}}{f}+\frac{\dot{V}^2}{V}-2\frac{\dot{f}\dot{V}}{f}=0
    \end{aligned}
\end{equation}
This second-order ODE system can be solved locally. However, the solutions to this ODE system will not yield a positively curved static triple with \( V = 0 \) on the boundary, as established by \textbf{Theorem \ref{thm_compact}}.\\

Finally, we prove \textbf{Theorem \ref{thm-vanishing}}
\begin{proof}
    For any element $[\omega]\in H^1(M)$, by the standard de Rham theory, we can find a 1-form $\omega\in [\omega]$ so that
    \begin{equation}
        \begin{aligned}
            d\omega &=\delta \omega=0 \text{ in } M\\
            \omega(\vec{n})&=0 \text{ on }\Sigma
        \end{aligned}
    \end{equation}
    where $\vec{n}$ is the outer normal vector. Let $X$ be the dual vector filed of $\omega$ w.r.t $g$, i.e. $\omega(\cdot)=g(X,\cdot)$. Then we have
    \begin{equation}
        \begin{aligned}\label{harmonic_X}
            &d\omega=0\longrightarrow g(\nabla_Y X,Z)=g(\nabla_ZX,Y)\quad  \forall Y,Z\\
            &\delta \omega=0\longrightarrow div(X)=0\\
           & \omega(\vec{n})=0\longrightarrow g(X,\vec{n})=0
        \end{aligned}
    \end{equation}
    We compute $\int_M\nabla^2V(X,X)dV_g$ as follows:
    \begin{equation}\label{integration_by_part}
        \begin{aligned}
            \int_M\nabla^2V(X,X)dV_g&=\int_Mg(\nabla_X\nabla V,X)\\
            &=\int_M div(g(\nabla V,X)X)-g(\nabla V,X)div(X)-g(\nabla V,\nabla_XX)\\
            &=\int_{\Sigma}g(\nabla V,X)g(X,\vec{n})-\int_M g(\nabla V,\nabla_XX)\\
            &=-\int_M g(\nabla V,\nabla_XX)
        \end{aligned}
    \end{equation}
    where we have used $div(X)=0$ in the third line and $g(X,\vec{n})=0$ in the fourth line. We continue by using another integration by part
    \begin{equation}
        \begin{aligned}
            \int_M\nabla^2V(X,X)dV_g&=-\int_{\Sigma}Vg(\nabla_XX,\vec{n})+\int_MVdiv(\nabla_XX)\\
            &=\int_MVdiv(\nabla_XX)
        \end{aligned}
    \end{equation}
    since $V=0$ on $\Sigma$. Pick local ON frame $\{e_i\}$ so that $\nabla e_i=0$ to compute $div(\nabla_XX)$.
    \begin{equation}\label{divXX}
        \begin{aligned}
            div(\nabla_XX)&=g(\nabla_{e_i}\nabla_XX,e_i)\\
            &=R(e_i,X,X,e_i)+g(\nabla_X\nabla_{e_i}X,e_i)+g(\nabla_{\nabla_{e_i} X-\nabla_{X} e_{i}}X,e_i)\\
            &=Ric(X,X)+g(\nabla_X\nabla_{e_i}X,e_i)+g(\nabla_{\nabla_{e_i} X-\nabla_{X} e_{i}}X,e_i)\\
            &=Ric(X,X)+I+II
        \end{aligned}
    \end{equation}
    The two terms on the RHS can be computed as follows:
    \begin{equation}\label{I}
            I=X(g(\nabla_{e_i}X,e_i))-g(\nabla_{e_i}X,\nabla_Xe_i)=0
    \end{equation}
    since $div(X)=0$ and $\nabla e_i=0$. And
    \begin{equation}
            II=g(\nabla_{\nabla_{e_i} X}X,e_i)=g(\nabla_{e_i} X,\nabla_{e_i} X)=|\nabla X|^2\label{II}
    \end{equation}
    We used first equation in (\ref{harmonic_X}). Now put (\ref{divXX})(\ref{I})(\ref{II}) into (\ref{integration_by_part}), we arrive at
    \begin{equation}
        \int_M\nabla^2V(X,X)=\int_M VRic(X,X)+|\nabla X|^2
    \end{equation}
    Use static equation (\ref{static_equ_2}), we get
    \begin{equation}
        0=\int_M-\Delta V|X|^2+|\nabla X|^2=\int_M nV|X|^2+|\nabla X|^2
    \end{equation}
    which implies $X=0$, and we are done.
\end{proof}

\textbf{Acknowledgements}

The author is thankful to Prof Xiaodong Wang for bringing this problem to my attention, Prof Lucas Ambrozio, Demetre Kazaras, Greg Galloway, and Eric Woolgar for helpful discussions, and Prof Stefano Borghini, Mattia Fogagnolo and Chao Xia for mentioning 
\cite{BF} to me.

\textbf{Data Availability Statement}
The authors declare that no datasets were generated or analyzed during the current study. This manuscript does not include any associated data.

\textbf{Conflict of Interest Statement}
I hereby declare that I have no financial, personal, or professional conflicts of interest related to this work.


\begin{thebibliography}{999}                                            
\bibitem{Am}Ambrozio, L. (2017). On static three-manifolds with positive scalar curvature. Journal of Differential Geometry, 107(1), 1-45.

\bibitem{Ander}Anderson, M.T. Static Vacuum Einstein Metrics on Bounded Domains. Ann. Henri Poincaré 16, 2265–2302 (2015). https://doi.org/10.1007/s00023-014-0367-8

\bibitem{BF}Borghini, S., \& Fogagnolo, M. (2023). Comparison geometry for substatic manifolds and a weighted Isoperimetric Inequality. arXiv preprint arXiv:2307.14618.

\bibitem{Bo}Bourguignon, J.P. (1975). Une stratification de l'espace des structures riemanniennes. Compositio Mathematica, 30, 1-41.

\bibitem{SB}Brendle, S. Constant mean curvature surfaces in warped product manifolds. Publ.math.IHES 117, 247–269 (2013). https://doi.org/10.1007/s10240-012-0047-5

\bibitem{BH}Brendle, S., \& Hung, P. K. (2024). Systolic inequalities and the Horowitz-Myers conjecture. arXiv preprint arXiv:2406.04283.

\bibitem{CG2}Cai, M., \& Galloway, G.J. (1999). Boundaries of zero scalar curvature in the AdS / CFT correspondence. Advances in Theoretical and Mathematical Physics, 3, 1769-1783.

\bibitem{CG}J. Cheeger, D. Gromoll, The splitting theorem for manifolds of nonnegative Ricci curvature.
J. Differ. Geom. 6, 119–128 (1971)


\bibitem{Co}J Corvino, Scalar curvature deformation and a glueing construction for the Einstein constraint equations,
Commun Math Phys 214 (2000) 137–189.

\bibitem{EH}Ellis GFR, Hawking SW. The Large Scale Structure of Space-Time. Cambridge University Press; 1973.

\bibitem{FF}Fédérer, H., \& Fleming, W.H. (1960). Normal and Integral Currents. Annals of Mathematics, 72, 458.

\bibitem{FM}A. Fischer and J. Marsden, Linearization stability of nonlinear partial differential equations, Proc. Symp. Pure Math. 27 (1975), 219-262.

\bibitem{Ga}Galloway, G.J. On the topology of black holes. Commun.Math. Phys. 151, 53–66 (1993). https://doi.org/10.1007/BF02096748

\bibitem{Ga1}Galloway, G. Maximum Principles for Null Hypersurfaces and Null Splitting Theorems. Ann. Henri Poincaré 1, 543–567 (2000). https://doi.org/10.1007/s000230050006

\bibitem{GSW}Galloway, G.J., Surya, S., \& Woolgar, E. (2002). On the Geometry and Mass of Static, Asymptotically AdS Spacetimes, and the Uniqueness of the AdS Soliton. Communications in Mathematical Physics, 241, 1-25.

\bibitem{Ko} O. Kobayashi, A differential equation arising from scalar curvature function, J. Math.
Soc. Japan 34, No. 4 (1982), 665-675.

\bibitem{La}J. Lafontaine, Sur la geometrie d’une generalisation de l’equation d’Obata, J. Math.
Pures Appliquees 62 (1983), 63-72.

\bibitem{Da}Lee, D. (2019). Geometric Relativity. Graduate Studies in Mathematics.

\bibitem{Li}Li, P. (2012). Geometric Analysis. Cambridge: Cambridge University Press.

\bibitem{LX}Li, J., \& Xia, C. (2019). An integral formula and its applications on sub-static manifolds. Journal of Differential Geometry, 113(3), 493-518.

\bibitem{Pe}Petersen, P. (2006). Riemannian geometry (Vol. 171, pp. xvi+-401). New York: Springer.

\bibitem{SY}Schoen, R., \& Yau, S. T. (1979). Existence of incompressible minimal surfaces and the topology of three dimensional manifolds with non-negative scalar curvature. Annals of Mathematics, 110(1), 127-142.

\bibitem{Wa}Wang, X. (2001). The Mass of Asymptotically Hyperbolic Manifolds. Journal of Differential Geometry, 57, 273-299.

\bibitem{Wa2}Wang, X. (2001). On Conformally Compact Einstein Manifolds. Mathematical Research Letters, 8, 671-688.

\bibitem{WY}Witten, E., \& Yau, S. (1999). Connectedness of the boundary in the AdS / CFT correspondence. Advances in Theoretical and Mathematical Physics, 3, 1635-1655.

\bibitem{Wo}Woolgar, E. (2016). The rigid Horowitz-Myers conjecture. Journal of High Energy Physics, 2017, 1-27.


\end{thebibliography}
\end{document}